\documentclass{article}
\usepackage{spconf,amsmath,graphicx,bm,amssymb}
\usepackage{algorithm}
\usepackage{color}
\usepackage{subfigure}
\usepackage{mathtools}

\usepackage[pdfstartview=FitH,
CJKbookmarks=true,
bookmarksnumbered=true,
bookmarksopen=true,
colorlinks,
pdfborder=001,
linkcolor=blue,
anchorcolor=blue,
citecolor=blue,
]{hyperref}
\hypersetup{hidelinks}

\usepackage{algorithmic}
\newcommand{\bI}{\mathbf{I}}

\newcommand{\bx}{\mathbf{x}}
\newcommand{\bX}{\mathbf{X}}

\newcommand{\by}{\mathbf{y}}
\newcommand{\bz}{\mathbf{z}}

\newcommand{\bD}{\mathbf{D}}

\newcommand{\bU}{\mathbf{U}}

\newtheorem{dingyi}{Definition~}%[section]
\newtheorem{dingli}{Theorem~}%[section]
\newtheorem{jiashe}{Assumption~}%[section]

\newcommand{\tnorm}[1]{\lVert\mkern-2mu |#1|\mkern-2mu\rVert}

\title{An efficient alternating Riemannian/projected gradient descent ascent algorithm
 for Fair Principal Component Analysis}
\name{Meng Xu$^{\star,\S}$, Bo Jiang$^{\dag}$, Wenqiang Pu$^{\ddag}$, Ya-Feng Liu$^{\star}$, and Anthony Man-Cho So$^{\sharp}$
}

\address{$^{\star}$LSEC, ICMSEC, AMSS, Chinese Academy of Sciences, Beijing, China \\
	$^{\S}$School of Mathematical Sciences, University of Chinese Academy of Sciences, Beijing, China\\
    $^{\dag}$Ministry of Education Key Laboratory of NSLSCS, \\
     School of Mathematical Sciences, Nanjing Normal University, Nanjing, China\\
    $^{\ddag}$Shenzhen Research Institute of Big Data, The Chinese University of Hong Kong, Shenzhen, China\\
    $^{\sharp}$Department of Sys. Eng. \& Eng. Mgmt, The Chinese University of Hong Kong, HKSAR, China\\
   Email:  xumeng22@mails.ucas.ac.cn, jiangbo@njnu.edu.cn, wenqiangpu@cuhk.edu.cn,\\ yafliu@lsec.cc.ac.cn, manchoso@se.cuhk.edu.hk 
}
  
\usepackage{cite}
\begin{document}

\ninept
\maketitle
%
%\wuhao {
\begin{abstract}
Fair principal component analysis (FPCA), a ubiquitous dimensionality reduction technique in signal processing and machine learning,  aims to find a low-dimensional representation for a high-dimensional dataset in view of fairness. The FPCA problem involves optimizing a non-convex and non-smooth function over the Stiefel manifold. The state-of-the-art methods for solving the problem are subgradient methods and semidefinite relaxation-based methods. However, these two types of methods have their obvious limitations and thus are only suitable for efficiently solving the FPCA problem in special scenarios. This paper aims at developing efficient algorithms for solving the FPCA problem in general, especially large-scale, settings. In this paper, we first transform FPCA into a smooth non-convex linear minimax optimization problem over the Stiefel manifold. To solve the above general problem, we propose an efficient alternating Riemannian/projected gradient descent ascent (ARPGDA) algorithm, which performs a Riemannian gradient descent step and an ordinary projected gradient ascent step at each iteration. We prove that ARPGDA can find an $\varepsilon$-stationary point of the above problem within $\mathcal{O}(\varepsilon^{-3})$ iterations. Simulation results show that, compared with the state-of-the-art methods, our proposed ARPGDA algorithm can achieve a better performance in terms of solution quality and speed for solving the FPCA problems. 

%
%[[[Will rewrite the abstract!!!]]]
\end{abstract}
\section{Introduction}
Dimensionality reduction, which aims to identify a low-dimensional subspace to represent high-dimensional data, has attracted much attention in recent decades due to its wide applications in signal processing and machine learning. 
The classical principal component analysis (PCA),  one of the most popular dimensional reduction techniques in practice, minimizes the single criterion of the total reconstruction error over the whole samples and may ignore the potential difference among different sensitive groups. 
Recently, the idea of fair machine learning has come to the stage, and applying this idea to PCA  gives rise to fair principal component analysis (FPCA) \cite{samadi2018price}.
It has been shown in \cite{samadi2018price,tantipongpipat2019multi,zalcberg2021fair} that compared to PCA, FPCA is more effective for not only machine learning applications such as credit scoring and image processing but also signal processing applications. 
%However, little of FPCA has been well understood, since FPCA can be formulated as a non-convex and non-smooth max-min problem which has been proved to be NP-hard \cite{tantipongpipat2019multi}. 
There is also growing interest in alternative definitions of fairness for PCA \cite{kamani2022efficient,lee2022fast,pelegrina2022analysis,pelegrina2022novel,vu2022distributionally}, 
%\cite{kamani2022efficient}-\cite{vu2022distributionally}, 
with various problem formulations and algorithms. 
%In this work, we focus on the widely used FPCA formulation \cite{samadi2018price}. 

However, FPCA is generally NP-hard (to solve exactly) \cite{tantipongpipat2019multi}, and there are two main approaches to the problem for approximate or local solutions. 
%One of the approach is the semidefinite relaxation (SDR) \cite{samadi2018price,tantipongpipat2019multi,olfat2019convex,morgenstern2019fair}, and some efficient algorithms have been designed. 
%For example, a multiplicative weight update is suggested in \cite{samadi2018price}.
%And in \cite{tantipongpipat2019multi}, an SDR extreme point, which can be obtained by an iterative rounding technique, is guaranteed to be a low-rank but still probably not an exact solution of FPCA  when the target dimension is much larger than the squared root of the number of groups. Although these algorithms are efficient for the fair machine learning setting, it is generally acknowledged that SDR scales poorly to high dimensions due to its large optimization variables needless to say the challenge to obtain a feasible solution. For FPCA in signal processing applications, with a totally different problem size, SDR is generally inapplicable in practice. To deal with such issues, the algorithms -based on subgradient come to the stage. For example, the Riemannian subgradient method in \cite{li2021weakly} can be used to directly solve FPCA. Meanwhile, the normal subgradient descent algorithm is suggested in \cite{zalcberg2021fair} to solve the unconstrained factorized FPCA, and any local minimizer of the factorized FPCA has been proved to be a global maximizer of FPCA under an orthogonality assumption. 
One approach is the semidefinite relaxation (SDR)-based algorithms \cite{samadi2018price,morgenstern2019fair,olfat2019convex,tantipongpipat2019multi}. In \cite{samadi2018price}, an SDR for FPCA was first proposed, and a multiplicative-weight style algorithm was used to approximately solve the semidefinite program. Some theoretical bounds for the SDR approach were also established in \cite{samadi2018price}.  Later on, the work \cite{tantipongpipat2019multi} improved the results in \cite{samadi2018price} to give a tighter approximation ratio of the SDR approach and showed that when the number of groups is two, the SDR approach can recover the solution of FPCA exactly. As mentioned in \cite{tantipongpipat2019multi}, solving the semidefinite program is the main bottleneck of the SDR approach. Hence, the SDR approach might be inapplicable to high-dimensional datasets \cite{zalcberg2021fair}.  The other approach to solving the FPCA problem is the subgradient methods. In particular, the work \cite{zalcberg2021fair} proposed an interesting equivalent unconstrained factorized FPCA formulation and applied the subgradient algorithm to solve the unconstrained equivalent problem. Notably, it has been shown in \cite{zalcberg2021fair} that any local minimizer of the factorized FPCA is a global maximizer of FPCA if the dataset satisfies an orthogonality condition. Therefore, the global solution of the original problem can be obtained if the orthogonality condition holds true.  
%Moreover, the Riemannian subgradient algorithm proposed recently in \cite{li2021weakly} can also be employed to solve FPCA. However, the subgradient type of methods has a low convergence rate which needs $\mathcal{O}(\varepsilon^{-4})$ iterations to get an $ \varepsilon $-nearly stationary point. 
Furthermore, the recently proposed Riemannian subgradient algorithm in \cite{li2021weakly} can be utilized for solving the FPCA problem. However, the Riemannian subgradient algorithm often exhibits a slow convergence rate, requiring approximately $\mathcal{O}(\varepsilon^{-4})$ iterations to return an $\varepsilon$-stationary point.

Motivated by the features of the existing approaches, it is of interest to develop more efficient as well as easy-to-implement methods for large-scale FPCA. By observing that FPCA can be rewritten as a non-convex linear minimax problem over the Stiefel manifold, in this paper, we propose an efficient alternating Riemannian/projected gradient descent ascent  (ARPGDA) algorithm for solving the FPCA problem and the more general non-convex linear minimax problem over  a Riemannian manifold.   
%In this paper, we focus on FPCA arising from signal processing, especially the filter design, as well as machine learning. 
%We first transform the non-smooth FPCA problem to a smooth minimax optimization over Stiefel manifold.
%We propose an efficient alternating Riemannian gradient (ARPGDA) algorithm, which
Our proposed ARPGDA algorithm only needs to perform a simple Riemannian gradient descent step and a projected gradient ascent step at each iteration, making it easy to implement and scalable to large-scale problems. 
In addition, we prove that our proposed ARPGDA algorithm has an iteration complexity of $ \mathcal{O}(\varepsilon^{-3}) $ for finding an $ \varepsilon$-stationary point (defined in Definition \ref{stationarypoint}), which matches the best known complexity for its Euclidean counterpart. 
%Additionally, to the best of our knowledge, our ARPGDA algorithm also fills the gap that no method with theoretical results has been designed for non-convex concave minimax optimization over Riemannian manifold before. Finally, simulation results show that our ARPGDA algorithm is more efficient than the algorithms based on the subgradient and SDR for filter design. Besides, the ARPGDA algorithm is also shown to be effective and promising for the applications in machine learning.
%To the best of our knowledge, our proposed ARPGDA algorithm is the first algorithm with an iteration complexity guarantee for general non-convex linear minimax optimization over Riemannian manifolds.  
Finally, simulation results show that, when applied to solve the FPCA problems from signal processing and machine learning applications, our proposed ARPGDA algorithm is more efficient in terms of both solution quality and speed compared to the subgradient and SDR methods.
% for filter design arising from signal processing and a block FPCA problem in machine learning. 

%It is worth mentioning some recent related works and comparing the algorithms proposed therein with our proposed ARPGDA algorithm. 
Our proposed ARPGDA algorithm can be seen as an extension of the alternating gradient projection algorithm \cite{pan2021efficient} for solving the non-convex linear minimax problem from the Euclidean space to the Riemannian manifold setting. 
To the best of our knowledge, the proposed algorithm is the first one with an iteration complexity guarantee for general non-convex linear minimax optimization over a Riemannian manifold.  Existing works on minimax optimization over a Riemannian manifold either focus on different problem settings or make strong assumptions that do not apply to the problems of our interest, including FPCA. Specifically, the works \cite{huang2023gradient} and \cite{jordan2022first,zhang2022minimax} considered the settings of non-convex strongly concave and geodesically convex geodesically concave minimax problems over Riemannian manifolds, respectively.
However, their settings do not apply to the FPCA problem. 
%Very recently, a Riemannian Hamiltonian gradient method was proposed in \cite{Andi2023Riemannian} to solve the Riemannian manifold constrained minimax problem, wherein a strong  Riemannian Polyak-\L ojasiewicz condition  on the Hamiltonian function need.  Nevertheless, the Riemannian Polyak-\L ojasiewicz condition is unlikely to hold in the FPCA problem.
Most recently, the work \cite{Andi2023Riemannian} introduced a Riemannian Hamiltonian gradient method for the Riemannian manifold constrained minimax problem, but it relies on a strong Riemannian Polyak-\L ojasiewicz condition that is unlikely to hold for the FPCA problem. 
%In contrast, our ARPGDA algorithm does not require any assumptions beyond general non-convexity, making it applicable to a broader range of problems such as FPCA.

% \rrev{???}
%It is also worthwhile mentioning some related works on minimax optimization. For the Euclidean space, an alternating \rev{projected gradient} algorithm for solving the non-convex concave minimax problem was proposed in \cite{xu2020unified}, and our proposed ARPGDA algorithm can be seen as an extension of it from the Euclidean space to the Riemannian manifold setting.
%Another closely related work is \cite{huang2023gradient}, which considered the non-convex strongly concave minimax problem over the Riemannian manifold. However, there is still no theoretical guarantee of their methods to solve our interested FPCA problem. Besides, very recently, \cite{zhang2022minimax,jordan2022first} also studied the general minimax problem over the manifold, however, they required that the objective function is geodesically convex geodesically concave, which does not apply to the FPCA problem.

%\vspace{-0.2cm}
\section{FPCA and Minimax Reformulation}
Let the columns of a $d\times N $ matrix $\bX = [\bx_1, \bx_2, \ldots, \bx_N]$ denote $ N $ data samples, each with $ d $ attributes. Suppose that the $N$ data samples belong to $ n $ groups, according to demographics or some other semantically meaningful clustering, and each group $i$ corresponds to a $d \times n_i$ submatrix $\bX_i$ with $\sum_{i = 1}^n n_i = N$.  
%For PCA, the goal is to find a subspace $ \bU\in\mathbb{R}^{d\times r} $, where $ r<d $, minimizes the total reconstruction error. 
Classical PCA aims to find a subspace with dimension $r < d$ to minimize the total reconstruction error or maximize the total variance equivalently. 
%However, some of the groups may suffer higher reconstruction errors than the others, which in turn may lead to bias for certain sensitive groups.
However, this may cause certain sensitive groups to suffer a higher reconstruction error than the others \cite{samadi2018price}.   
To reduce such disparity, FPCA minimizes the maximum reconstruction error among the $ n $ groups, which is equivalent to maximizing the minimum variance among the $n$ groups. 
%which is equivalent to maximize the minimum of the Frobenius norms, $ \| \bU\bU^T\bX_i\|^2_F$. 
Formally, FPCA can be formulated as \cite{tantipongpipat2019multi}
\begin{equation}\label{FPCAblock}
\text{FPCA:}\quad\max_{\bU\in\mathcal{S}}\min_{i=1,2,\ldots,n}\ f_i(\bU):=\langle \bX_i\bX_i^T, \bU\bU^T\rangle,
\end{equation}
where $\langle \cdot, \cdot \rangle$ is the standard Euclidean inner product %on $\mathbb{R}^{d \times \rev{d}}$
 and $\mathcal{S} := \{\bU \in \mathbb{R}^{d \times r}\mid\bU^T \bU = \bI_r\}$ is known as the Stiefel manifold. 
%When the number of samples $ N $ is much larger than $ n $, the matrix $ \bX $ is divided into several blocks. Therefore, we also call this problem block-FPCA, which is introduced in the fair machine learning \cite{tantipongpipat2019multi}. We note that our FPCA is different from one in \cite{samadi2018price}, but equivalent by additive and multiplicative scalings.
Note that if $n = 1$, problem (\ref{FPCAblock}) reduces to classical PCA. If $n\ll N$, %the number of samples $N$ is much larger than $n$, 
problem (\ref{FPCAblock}) is called block FPCA \cite{zalcberg2021fair}, which was first introduced in the context of fair machine learning \cite{samadi2018price}.  If $n = N$, problem (\ref{FPCAblock}) becomes 
%As detailed in the introduction, FPCA also arises in signal processing applications, such as filter design \cite{zalcberg2021fair}, where the FPCA problems have totally different shapes. 
%%Consider a problem in filter design for detection. 
%Let $\{\bx_i\}_{i=1}^n \subset \mathbb{R}^d$ denote the high resolution digital signals 
%%(or even an infinite dimension analog signal) 
%of $ n $ kinds of materials. The aim is to design a few expensive sensors with low dimensions that can downsample the high dimensional signals effectively. Therefore, it is practical to maximize the weakest reflectance, i.e. 
\begin{equation}\label{FPCA}
	\max_{\bU\in\mathcal{S}}\min_{i=1,2,\ldots,n}\ \langle \bx_i\bx_i^T, \bU\bU^T\rangle,
\end{equation}
which arises in many signal processing applications,  such as multicast downlink beamforming \cite{sidiropoulos2006transmit,luo2010semidefinite,lu2017anefficient, chen2018joint} and filter design \cite{zalcberg2021fair}. In the context of filter design, $\{\bx_i\}_{i=1}^n$ denote the high-resolution digital signals of $n$ kinds of materials, 
and the parameters always satisfy $r < n = N \ll d$, which is sharply different from the settings in fair machine learning, in which $d$ might be smaller than $N$.
Let $\Delta:=\left\{\by\in\mathbb{R}^n\ |\ \sum_{i=1}^{n} y_i =1,\ y_i\geq0,\ i=1,2,\ldots,n\right\}$ be the standard simplex.  The FPCA problem (\ref{FPCAblock}) can be equivalently rewritten as follows:
\begin{equation}\label{maxminreform}
\min_{\bU\in\mathcal{S}}\max_{\by\in\Delta}\ \sum_{i=1}^{n}y_i (-f_i(\bU)).
\end{equation}  
 %More generally, the above problem can be considered as a non-convex concave minimax problem on the Riemannian manifold, i.e.
Hence, we generalize (\ref{maxminreform}) to consider the non-convex linear minimax problem over a Riemannian submanifold $\mathcal{M}$ of a finite-dimensional Euclidean space $\mathcal{E}_1$, i.e., 
\begin{equation}\label{minimaxRM}
	\min_{\bU\in\mathcal{M}} \max_{\by\in\mathcal{Y}}\ f(\bU,\by),
\end{equation}
where $ f(\bU,\by):\mathcal{M}\times\mathcal{Y}\to \mathbb{R} $ is a smooth function that is linear with respect to $ \by $ but possibly non-convex with respect to $ \bU $, and $ \mathcal{Y} $ is a non-empty compact convex set in a finite-dimensional Euclidean space $\mathcal{E}_2$. In the rest of this paper, we shall develop an efficient algorithm for solving the general non-convex linear minimax problem (\ref{minimaxRM}) with the Riemannian manifold constraint, which includes our interested FPCA problem as a special case.

\section{Proposed ARPGDA Algorithm}
In this section, we propose an efficient and scalable ARPGDA algorithm for solving the non-convex linear minimax problem (\ref{minimaxRM}) and establish its iteration complexity. 

\subsection{Proposed Algorithm}
% More specifically, we perform a simple Riemannian gradient descent step to update $\bU$ and an ordinary \rev{projected gradient}  ascent step to update $\by$. 
%Besides, the high dimension of the matrix variable $ \bU $ in the real world also poses a great challenge to the algorithm. 
%To handle with the T scale problem, we adopt an alternating first order single-loop approach to update $ \bU $ and $ \by $, which has a low computational complexity in each iteration. 
%In order to tackle the non-convex constraint, we note that the Riemannian manifold is a special non-convex feasible set with relatively benign properties, and both the theory and algorithm of Riemannian manifold optimization are well studied. 
%Therefore, we apply the simple Riemannian gradient descent method when updating $ \bU $ and apply the \rev{projected gradient} method when updating $ \by $.

%Our proposed algorithm can be regarded as an extension of the algorithms proposed in \cite{xu2020unified} and \cite{HiBSA}, which are designed for  minimax problems with convex constraints and thus can not be applied directly to our interested problems. Similar to \cite{xu2020unified} and \cite{HiBSA}, we introduce a regularized 
Before presenting the ARPGDA algorithm, we introduce some basic objects associated with the Riemannian manifold $\mathcal{M}$. Let $ T_{\bU}\mathcal{M} $ denote the tangent space to $\mathcal{M}$ at $\bU\in\mathcal{M}$. The Riemannian manifold $ \mathcal{M} $ is endowed with a smooth inner product $ \langle\cdot,\cdot\rangle_{\bU}\colon T_{\bU}\mathcal{M}\times T_{\bU}\mathcal{M}\to \mathbb{R} $. 
Let  $\mathcal{R}_{\bU}\colon T_{\bU} \mathcal{M} \to \mathcal{M}$ be the retraction at $\bU$, which is a smooth map satisfying i) $\mathcal{R}_{\bU}(\mathbf{0}_{\bU}) = \bU$, where $\mathbf{0}_{\bU}$ is the zero element in $T_{\mathbf{U}} \mathcal{M}$; ii) the differential of $\mathcal{R}_{\bU}(\mathbf{0}_{\bU})$ is the identity map. 
%\rrev{\sout{A classical retraction is the  exponential mapping.}}
For a smooth function $f\colon \mathcal{M} \to \mathbb{R}$, the Riemannain gradient $\mathrm{grad} f(\bU)$ is the unique tangent vector at $\bU$ satisfying $\langle \mathrm{grad} f(\bU), \bD\rangle_{\bU} =\mathrm{D} f(\bU)[\bD]$, where $\mathrm{D} f(\bU)[\bD]$ denotes the directional derivative of $f$ along the direction $\bD$. For the Stiefel manifold,  the endowed norm $\|\cdot\|_{\bU}$ is taken as the usual Euclidean norm $\|\cdot\|$ and the Riemannian gradient is
\begin{equation}\label{Riemanniangradient}
	\mathrm{grad} f(\bU) = \nabla f(\bU) - \bU (\bU^T \nabla f(\bU) + \nabla f(\bU)^T \bU)/2.
\end{equation}
Moreover, we adopt the polar decomposition-based retraction \cite{absil2008optimization}, which satisfies $\|\mathcal{R}_{\bU}(\bD) - \mathcal{R}_{\bU}(\mathbf{0}_\bU)\| \leq  \|\bD\|$ and $\|\mathcal{R}_{\bU}(\bD) - \mathcal{R}_{\bU}(\mathbf{0}_\bU) - \bD \|\leq 1/2\|\bD\|^2$;  see \cite[Lemma A.2]{jiang2017vector} for details. 
 
Recall that (\ref{minimaxRM}) is a non-convex linear minimax problem with a non-convex function in $\bU$ and a manifold constraint $\mathcal{M}$, which is difficult to solve.
To handle the large-scale problem, we adopt an alternating first-order approach to update $\bU$ and $\by$, which has a very low cost at each iteration. 
%Some efficient first-order methods for solving general non-convex concave minimax problems were recently proposed in \cite{xu2020unified} and \cite{HiBSA}. However, the methods therein cannot be directly applied to solve our problem (4) due to the manifold constraints. 
Inspired by the regularization techniques in %\cite{pan2021efficient}
\cite{nouiehed2019solving,lu2020hybrid,pan2021efficient,xu2023unified}, we introduce a regularized version of the original objective function in (\ref{minimaxRM}) at the $k$-th iteration, i.e.,
\[
f_k(\bU,\by):=f(\bU,\by)-\frac{\lambda_k}{2}\left\|\by\right\|^2, 
\]
%\begin{equation*}
%\tilde{f}_k(\bU,\by)=f(\bU,\by)-\frac{\lambda_k}{2}\left\|\by\right\|^2-\frac{\beta_k}{2}\|\by-\by_k\|^2,
%\end{equation*} 
where $ \lambda_k\geq0 $ is a regularization parameter.
%and $\|\cdot\|$ denotes the standard norm in the Euclidean space $\mathcal{E}_2$.
%\rrev{\sout{It is worth mentioning that there is a vast literature} \cite{nouiehed2019solving,lu2020hybrid,xu2023unified} \sout{highlighting the necessity of the regularized term for the convergence of the gradient descent/ascent type of algorithms for solving even convex concave  minimax problems.}} \rrev{Comment by Bo: For convex concave minimax problem, there is no need to add the regularized term, therefore, it is better to delete the sentence.}
%hommes2012multiple
Our proposed algorithm is based on the above regularized function $f_k$ at the $k$-th iteration. In particular, the $k$-th iteration of the proposed algorithm performs a Riemannian gradient descent step to update $ \bU $, followed by a projected gradient ascent step to update $ \by $.
%Since the above algorithm performs Riemannian gradient descent steps and  \rev{projection} gradient ascent steps in an alternating fashion, 
Hence,  we name it the alternating Riemannian/projected gradient descent ascent (ARPGDA) algorithm. 
 %More specifically, 
%\rrev{\sout{solves the inner maximization subproblem $\max_{\by\in\mathcal{Y}} \tilde{f}_k(\bU,\by)$ to update $\by$. Note that the inner maximization is a strongly-concave quadratic programming, which can be exactly solved by a \rev{projected gradient} step, i.e.,
%More specifically, for a given pair $ (\bU_k,\by_k) $, we perform a Riemannian gradient descent step with respect to $\bU$. 
%to get a temporary variable $ \bD_k $. And then calculate the retraction of $ \bD_k $ as $ \bU_{k+1} $.
%We use the Forbenius norm for a matrix and the $ l_2-$norm for a vector.
%Then our proposed algorithm uses a Riemannian gradient descent step of $ \Phi_{k} $ to update $\bU$.}} 
Let $\mathrm{grad}_{\bU}f(\bU,\by)$ be the Riemannian gradient with  respect to $\bU$, which can be calculated according to \eqref{Riemanniangradient}. Noting that $\mathrm{grad}_{\bU} f_{k}(\bU_k, \by_{k}) = \mathrm{grad}_{\bU} f\left(\bU_{k}, \by_{k}\right) $,
we update $\bU_{k+1}$ as follows:  
%\begin{equation*}
%	\begin{aligned}
%		\bD_{k}=&\arg \min _{\bD \in T_{\bU_k}\mathcal{M}}\left\langle\mathrm{grad}_{\bU} \tilde{f}_k\left(\bU_{k}, \by_{k}\right), \bD\right\rangle_{\bU_k}+\frac{\beta_{k}}{2}\left\|\bD\right\|^{2}_{\bU_k}\\
%		=&-\frac{1}{\beta_k}\mathrm{grad}_{\bU} f\left(\bU_{k}, \by_{k}\right),\\
%		\bU_{k+1}=&\ \mathcal{R}_{\bU_k}(\bD_k),
%	\end{aligned}
%\end{equation*}
\begin{equation}\label{updateU}
		\bU_{k+1}=\mathcal{R}_{\bU_k}( -\zeta_k\mathrm{grad}_{\bU} f\left(\bU_{k}, \by_{k}\right)),
\end{equation}
%where \rrev{\sout{$ \bD_{k} = -\zeta_k\mathrm{grad}_{\bU} f\left(\bU_{k}, \by_{k+1}\right) $ and}} 
where $ \zeta_k> 0  $ is the stepsize. 
With the identity $\nabla_{\by} f_k(\bU_{k+1}, \by_k) = \nabla_{\by} f\left(\bU_{k+1}, \by_{k}\right)- \lambda_k \by_{k}$, we  update $\by_{k+1}$ as follows:
\begin{equation}\label{updatey}
	\by_{k+1}=\mathcal{P}_{\mathcal{Y}}\Big(\by_{k}+\frac{1}{\lambda_k+\beta_k} \big(\nabla_{\by} f\left(\bU_{k+1}, \by_{k}\right)- \lambda_k \by_{k}\big)\Big),
\end{equation}
 where $ \mathcal{P}_{\mathcal{Y}}(\bz) $ denotes the projection of $\bz$ onto the set $\cal Y$ and $\frac{1}{\lambda_k+\beta_k}> 0$ is the stepsize.
%\rrev{\sout{parameter of updating $ \bU $}}
%\rrev{\sout{For our interested Stiefel manifold, we can calculate $ \bD_k $ via (\ref{Riemanniangradient}), and the QR decomposition is commonly used to compute the retraction, i.e., $ \mathcal{R}_{\bU}(\bD) $ is the $ \bQ $ matrix in $ \bD=\bQ\bR $.}}
The ARPGDA algorithm is formally presented in Algorithm \ref{algonp}.

\vspace{-0.1cm}
\begin{algorithm}
	\caption{ARPGDA Algorithm for Solving Problem (\ref{minimaxRM})}
	\begin{algorithmic}\label{algonp}
		\small
		\STATE Step 1 \ Input $\bU_1, \by_1$,  $\varepsilon, \{\lambda_k\}, \{\beta_k\}, \{\zeta_k\}$; set $k=1$.
		\STATE Step 2 \ Calculate $ \bU_{k+1} $ via (\ref{updateU}).
		\STATE Step 3 \ Calculate $ \by_{k+1} $ via (\ref{updatey}).
		\STATE Step 4 \ If $\mathcal{E}(\bU_{k+1},\by_{k+1})\leq \varepsilon$ (see \eqref{FNE2} further ahead), stop; otherwise, set $k\gets k+1$ and go to Step 2.
	\end{algorithmic}
\end{algorithm}

\vspace{-0.2cm}

%Observe that the proposed ARPGDA algorithm is a single-loop first order algorithm, and it directly uses $ \bU $ to calculate steps which protects the low dimension property in our interested FPCA problems. Therefore, the obtained solution is guaranteed to be feasible, and Algorithm \ref{algonp} has less computational complexity than the SDR approach which introduces a high-dimensional matrix variable $ \bU\bU^T $. In fact, the update of $ \by $ is also efficient in our interested FPCA problems since a fast algorithm for projection onto the simplex, with the complexity of $ \mathcal{O}(n) $ in practice, has already been proposed \cite{projection}. The total computational complexity of ARPGDA for each iteration is $ \mathcal{O}((2nr+2r^2)d) $ which is a bit higher than those of the subgradient type methods in \cite{zalcberg2021fair} and \cite{li2021weakly}, which are $ \mathcal{O}((nr+3r^2)d) $ and $ \mathcal{O}((nr+2r^2)d) $ respectively. 

Our proposed ARPGDA algorithm is a first-order algorithm, as it only uses the current gradient information of $f(\bU,\by)$ to generate the next point. Hence, it is suitable to solve our interested large-scale FPCA problem (\ref{FPCAblock}). When applied to FPCA, each iteration of the ARPGDA algorithm essentially performs a Riemannian gradient step to inexactly solve a classical PCA problem with an adaptive weighted matrix $ \sum_{i=1}^n y_i \bX_i \bX_i^T $ (where $y_i$ can be seen as the weight of the data in group $i$), 
%\rrev{\sout{In particular, when ARPGDA is applied to solve FPCA problem (\ref{FPCAblock})}}
and the per-iteration cost of updating $\bU$ and $\by$ is $\mathcal{O}(Nrd + r^2d)$ and $\mathcal{O}(n\mathrm{log}n)$ \cite{held1974validation,condat2016fast}, respectively.  Compared with the SDR approach for FPCA, ARPGDA has a much lower computational complexity since the SDR approach works with a high-dimensional matrix variable in $\mathbb{R}^{d \times d}$.
%and the cost is $O(Nrd + d^{6.5})$ \cite{tantipongpipat2019multi}. 
Compared with the subgradient-type methods \cite{zalcberg2021fair} and \cite{li2021weakly}, ARPGDA has a similar per-iteration cost.  However, as shown later in the simulation, we can see that ARPGDA always converges significantly faster than the subgradient methods. 

\subsection{Iteration Complexity}
%We start to discuss the theoretical properties of the proposed ARPGDA algorithm. Before presenting the convergence results, we need to propose an alternative optimality condition for minimax problem (\ref{minimaxRM}). In fact, the optimality for non-convex minimax problems has been an important research front since finding a Nash equilibrium, the  classical optimality, is NP-hard due to the non-convexity. 
Before presenting the iteration complexity result for the ARPGDA algorithm, we first define the notion of an $\varepsilon$-stationary point of  problem (\ref{minimaxRM}) as follows, which is motivated by the definition in \cite{nouiehed2019solving,pan2021efficient} for the convex constrained counterpart of problem (\ref{minimaxRM}).
\begin{dingyi}\label{stationarypoint}
	For any given $\varepsilon>0$, we say that $(\tilde \bU,\tilde \by)$ is an $\varepsilon$-stationary point of problem (\ref{minimaxRM}) if 
	%\begin{equation*}\label{FNE1}
	%	\|\mathrm{grad}_{\bU} f(\tilde \bU, \tilde \by)\|\leq\varepsilon,
	%\end{equation*} and
	% there exists some $\gamma > 0$ such that
	%\begin{equation*}\label{FNE2}
		$\mathcal{E}(\tilde \bU, \tilde \by)\leq\varepsilon,$
		where 
		\begin{equation}\label{FNE2}
		\mathcal{E}(\tilde \bU, \tilde \by):= \max\!\big\{\|\mathrm{grad}_{\bU} f(\tilde \bU, \tilde \by)\|_{\tilde \bU},  \max_{\by \in \mathcal{Y}} \langle \nabla_{\by} f(\tilde \bU, \tilde \by),\by - \tilde \by\rangle\big\}.
		\end{equation}
	%\end{equation*}
%	$
%		\nabla_{\bU} G(\bU,\by):=-\mathrm{grad}_{\bU} f(\bU,\by),
%		\nabla_\by G(\bU,\by):=\by - \mathcal{P}_\mathcal{Y}\left(\by_{}+ \nabla_{\by} f\left(\bU_{}, \by_{}\right)\right)$, and $
%		\left\| \nabla G(\bU,\by)\right\|^2 := \left\|\nabla_\bU G(\bU,\by)\right\|_\bU^2+\left\|\nabla_{\by} G(\bU,\by)\right\|^2$.
\end{dingyi}
%Similar to \cite{xu2020unified} and \cite{HiBSA}, we define the $ \varepsilon-$stationary point to be the measure of optimality.
%We define the $\epsilon$-stationary point of problem (\ref{minimaxRM}) as follows. 
%\textcolor{red}{Define}  
%\begin{equation}
%	\nabla \mathbb{G}(U,y): = 
%	\begin{bmatrix}
%		\mathrm{grad}_U f(U,y) \\ 
%		y - \mathcal{P}_Y(y + \nabla_y f(U,y))
%	\end{bmatrix}
%\end{equation}. 
%\begin{dingyi}\label{stationarypoint}
%	For any given $\varepsilon>0,$ we say that $(\bU^*,\by^*)$ is an $ \varepsilon$-stationary point of problem (\ref{minimaxRM}) if %the following holds
%	 $\left\| \nabla G(\bU^*,\by^*)\right\|\leq\varepsilon,$ where
%	 \begin{equation*}
%	 	\nabla_{\bU} G(\bU,\by):=-\mathrm{grad}_{\bU} f(\bU,\by),
%	 \end{equation*}
%	 \begin{equation*}
%	 	\nabla_\by G(\bU,\by):=\by - \mathcal{P}_\mathcal{Y}\left(\by_{}+ \nabla_{\by} f\left(\bU_{}, \by_{}\right)\right),
%	 \end{equation*}
%	 \begin{equation*}
%	 	\left\| \nabla G(\bU,\by)\right\|^2 := \left\|\nabla_\bU G(\bU,\by)\right\|_\bU^2+\left\|\nabla_{\by} G(\bU,\by)\right\|^2.
%	 \end{equation*}
%\end{dingyi}
%If $ \left\| \nabla G(\bU,\by) \right\|=0 $, we can conclude that $ \left\|\left(\nabla G(\bU,\by)\right)_{\bU}\right\| $ and $ \left\|\left(\nabla G(\bU,\by)\right)_{\by}\right\| $ are zeros, which are the optimality conditions for the Riemannian optimization with respect to $ \bU $ and the convex optimization with respect to $ \by $. It has been proved that $ \left\| \nabla G(\bU,\by) \right\|=0 $ is also equivalent to the KKT conditions of problem (\ref{maxminreform}) \cite{HiBSA}.

We also need to make the following (somewhat standard) assumption about the smoothness of the function $ f $.
\begin{jiashe}\label{A1}
	The function $f$ is continuously differentiable, and there exist positive constants $L_1~\text{and}~L_2$ such that for any 
	 $ \bU\in\mathcal{M}, \bD\in T_{\bU}\mathcal{M}$,
	 and $\bar{\bU}=\mathcal{R}_{\bU}(\bD)$, we have
	\begin{subequations}
	%&\left\|\operatorname{grad}_{\bU} f\left(\bU_{1}, \by\right)- T_{\bar{\bU}}^{\bU}\operatorname{grad}_{\bU} f\left(\bar{\bU}, \by\right)\right\| \leq L_1\|\bD\|,\\ 
	\begin{equation}\label{L11}
		f(\bar{\bU},\by)\leq f(\bU,\by)+\left\langle\mathrm{grad}_{\bU} f(\bU,\by),\bD\right\rangle_{\bU}+\frac{L_1}{2}\left\|\bD\right\|_{\bU}^2,
	\end{equation}
	\begin{equation}\label{L12}
		\left\|\operatorname{grad}_{\bU} f\left(\bU, \by\right)-\operatorname{grad}_{\bU} f\left(\bU, \bar{\by}\right)\right\|_{\bU} \leq L_2\left\|\by-\bar{\by}\right\|.
	\end{equation}
%	\begin{equation}
%		\tag{b}\left\|\nabla_{\by} f\left(\bU, \by\right)-\nabla_{\by} f\left(\bar{\bU}, \by\right)\right\| \leq \rev{L_2}\|\bD\|_{\bU},
%	\end{equation}
%	\begin{equation}
%		\tag{c}\left\|\nabla_{\by} f\left(\bU, \by\right)-\nabla_{\by} f\left(\bU, \bar{\by}\right)\right\| \leq L_2\left\|\by-\bar{\by}\right\|,
%	\end{equation}
\end{subequations}
%Moreover,  we asume that for any $ \bU\in\mathcal{M}$ with $\bD\in T_{\bU}\mathcal{M}$, there exist a constant $L_1$ such that
\end{jiashe}
%Assumption \ref{A1} is the gradient $L$-continuous condition for Riemannian optimization which is commomly used in minimax optimization \cite{xu2020unified,lin2020gradient} and Riemannian optimization \cite{huang2023gradient,sato2019riemannian}. And it implies that the partial Riemannian gradient
%$ \mathrm{grad}_\bU f(\bU,\by) $ for all $ \by \in \mathcal{Y} $ is retraction $ L_1 $-smooth as in \cite{huang2023gradient,sato2019riemannian}, which provides a sufficent descent condition.
Condition \eqref{L11} is  a standard assumption in manifold optimization \cite{boumal2019global}, and condition \eqref{L12} is widely used in the minimax optimization literature \cite{pan2021efficient,huang2023gradient,xu2023unified}.  
For our interested FPCA problem \eqref{maxminreform}, by some calculations, we can show that Assumption \ref{A1} holds with $L_1 = 2\max_i \|\bX_i\bX_i^T\|$ and $L_2 = 2\sqrt{\tnorm{\sum_{i=1}^n \bX_i\bX_i^T \bX_i\bX_i^T}_{(r)}}$. Here, $\tnorm{\cdot}_{(r)}$ is the Ky Fan $r$-norm, which is the sum of the $r$ largest singular values of a matrix. %\rrev{Bo: Meng, please check the two constants.}\rrev{\sout{see \cite{jiang2017vector} for details.}}

Now we are ready to present the iteration complexity result of the ARPGDA algorithm for solving a general non-convex linear minimax problem over a Riemannian manifold, namely, problem \eqref{minimaxRM}. 
%Since whether the objective function $f$ is strongly concave with respect to $\by$ is crucial to the convergence performance,
%we will discuss the two cases, namely, the strongly concave and concave cases, separately. 
%Let $ T(\varepsilon):=\min\left\{ k \ |\  \left \| \nabla G(\bU_k,\by_k) \right \| \leq \varepsilon \right \} $ denote the minimum number of iterations to get an $ \varepsilon$-stationary point. We have the following theorems.
\begin{dingli}\label{TheoremLinear}
	Suppose that Assumption \ref{A1} holds and the parameters in Algorithm \ref{algonp} satisfy $\lambda_k= \frac{\varepsilon}{8R^2}$, $\beta_k=\mu_k k^{-\rho}, 0 \leq \mu_k\leq\mu, \zeta_k = \theta/(L_1+\frac{L_2^2}{\lambda_k+\beta_k + \beta_{k+1}})$, where $\rho>1, \theta \in (0,2)$ are constants, and $R=\max_{\by\in\mathcal{Y}}\|\by\|$.
	%where $R=\max_{\by\in\mathcal{Y}}\|\by\|, \rho>1$.
	%Then there exists  $ (\frac{8R^2\mu}{\varepsilon})^{\frac{1}{\rho}}<\hat{k}\leq T $, where  such that $ (\bU_{\hat{k}},\by_{\hat{k}}) $ is an $\varepsilon$-FNE of problem (\ref{minimaxRM}). 
	Then, Algorithm \ref{algonp}  is guaranteed to return an $\varepsilon$-stationary point of problem (4) within 
	%can return an $\varepsilon$-stationary point of problem (\ref{minimaxRM}) 
	 $\mathcal{O}(\varepsilon^{-3}) $ iterartions. 
\end{dingli}
Our  iteration complexity result of $\mathcal{O}(\varepsilon^{-3})$ in Theorem \ref{TheoremLinear} matches the best known complexity result of general non-convex linear minimax optimization, in which the Riemannian manifold is replaced by a non-empty compact convex set \cite{pan2021efficient}. Due to the space limitation, instead of giving a rigorous proof of Theorem \ref{TheoremLinear}, we just give some key steps of the proof.
%we just sketch the proof of the theorem. 
%\rrev{\sout{We prove that $\bU_{k}$ and $ \by_{k+1} $ generated by Algorithm \ref{algonp} will satisfy the conditions in Definition \ref{stationarypoint} respectively.}}
%The key insight is that the function $\Phi_{\lambda\beta}(\bU)$ can be  considered as a smoothed function of $ \tilde{\Phi} $ in (\ref{minimaxRM}).
%Since the inner maximization is exactly solved, the key is to analyze the convergence of $\{\bU_{k}\}_{k\geq1} $. 
First, based on the update schemes  \eqref{updateU} and \eqref{updatey}, we can estimate the sufficient reduction of $f_{k}(\bU,\by)$ at each iteration as follows:
\begin{align*}
	& f_{k+1}(\bU_{k+1}, \by_{k+1}) - f_k(\bU_{k}, \by_{k})\\
	\leq&-\frac{2 - \theta}{2\theta}\zeta_k\|\mathrm{grad}_\bU f(\bU_k,\by_{k})\|_{\bU_k}^2+ \frac12(\lambda_k-\lambda_{k+1}+4\beta_k)R \\
	& -\frac12(\beta_{k} \|\by_{k} -\by_{k-1}\|^2 - \beta_{k+1}\|\by_{k+1} -\by_{k}\|^2).
\end{align*}
 %\rrev{\sout{and show there exists some $\bU_{\hat{k}}$ satisfies (\ref{FNE1}) in $ \mathcal{O}(\varepsilon^{-3}) $ iterations.}}
Second, due to the fact that $f$ is linear in  $\by$, we have
\[
\max_{\by \in \mathcal{Y}} \langle \nabla_{\by} f(\bU_{k}, \by_{k}),\by - \by_{k}\rangle \leq 4R^2(\lambda_k + \beta_k).
\]
 %$(\lambda_k + \beta_k)\|\by_{k+1} - \mathcal{P}_\mathcal{Y}(\by_{k+1}+ \frac{1}{\lambda_k + \beta_k}\nabla_{\by} f(\bU_k, \by_{k+1})\| \leq (\lambda_k + 2\beta_k) R$.
Finally, by carefully choosing the parameters, we can complete the proof of Theorem  \ref{TheoremLinear}.
%First, we estimate a potential function such that it has a sufficient reduction at each iteration. 
It should be mentioned that our key proof process is motivated by and similar to that in \cite{pan2021efficient}. The main difference lies in the first step, where we use conditions \eqref{L11} and \eqref{L12} to bound the change of the function value. 
Even in the Euclidean setting, our condition \eqref{L11} is still weaker than the first Lipschitz continuous gradient assumption of Assumption 1 in \cite{pan2021efficient}. 
\section{Simulation Results}
%In this section, we illustrate both the effectiveness and efficiency of the proposed ARPGDA algorithm with numerical experiments. We test the algorithms on our interested FPCA problems from signal processing and machine learning applications.
In this section, we present numerical results  %from signal processing and machine learning applications
to show both the effectiveness and efficiency of the proposed ARPGDA algorithm when solving the FPCA problem. We mainly compare ARPGDA with the subgradient (SG) method, which is designed to solve the unconstrained factorized FPCA problem \cite{zalcberg2021fair}, and the Riemannian subgradient (RSG) method \cite{li2021weakly}, which solves the following equivalent formulation of problem \eqref{FPCAblock}:
 \[
 \max_{\bU\in\mathcal{S}}\ \Phi(\bU) := \min_{i = 1,2, \ldots, n} f_i(\bU).
\] The stepsize of the two subgradient methods is $c/k^{1/2}$, where $c$ is chosen carefully in each test. We stop SG or RSG when $\Phi(\bU) \geq  (1 - 10^{-4}) \Phi_{\text{ARPGDA}}$ or the iteration number hits the maximum number of iterations $10^5$. Here, $\Phi(\bU) $ and $\Phi_{\text{ARPGDA}}$ denote the objective values returned by SG or RSG  and ARPGDA, respectively. When the dimension $ d $ of the problem is not high, we also compare the SDR-based approach \cite{tantipongpipat2019multi}. All results in this section are averaged over 10 runs with different randomly generated initial points.
%as $\max_{\mathbf{\bU} \in \mathcal{S}} \Phi(\bU)$, with $\Phi(\bU)$ defined as (\ref{Phi}). %with $\Phi(\bU) = \min_{i = 1,\ldots n} f_i(\mathbf{U})$. 
 %We do not present the results of the  because it is inapplicable to this high-dimensional setting.

\textbf{Numerical Results on FPCA Problem (1) with $n=N$.}
The FPCA problem (\ref{FPCAblock}) with $n=N$, namely, problem \eqref{FPCA},  arises in many signal processing applications \cite{sidiropoulos2006transmit,zalcberg2021fair}.  For ARPGDA, we 
%to solve the reformulation (\ref{maxminreform}) with 
choose $\varepsilon = 10^{-3} \max_{i} \|\bx_i\|^2$,  $\rho=1.1$, $\theta = 1.5$,   %\rrev{\sout{$\mu_k$ to a piecewise function, i.e., when $k\leq500$,}} 
and $\mu_k \equiv 30n^2\sqrt{r}$. Note that for FPCA, the constant $R = 1$.
% if $k\leq500$ 
%and set $\mu_k=3n^2$ otherwise. 

\begin{figure}[t]
	\centering
	\begin{subfigure}
		\centering
		\includegraphics[scale=0.52]{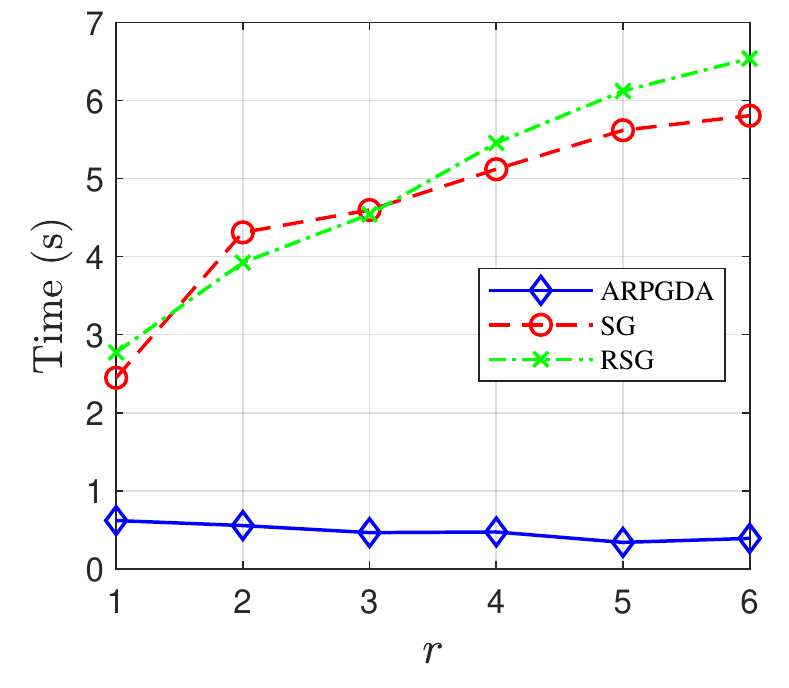}
		\vspace{-7pt}
	\end{subfigure}
	\begin{subfigure}
		\centering
		\includegraphics[scale=0.52]{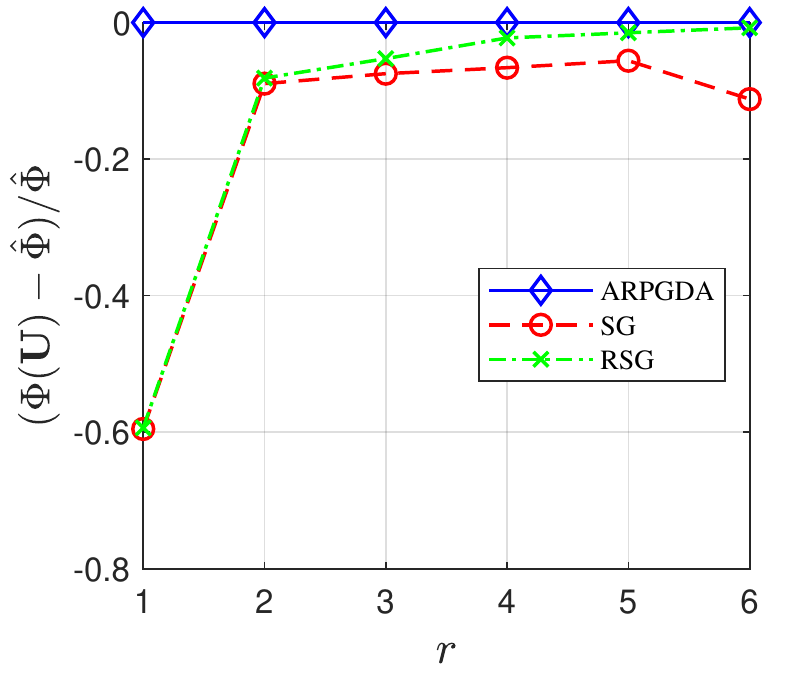}
		\vspace{-7pt}
	\end{subfigure}
	\caption{Average performance comparison on the synthetic data. %The right figure presents the quality of the obtained solutions. %Since the approximate solutions of SG  are not feasible, the projection to feasible region and the low efficiency may be the reason of the worse objective function values.
	}
	\vspace{-7pt}
	\label{Synpic1}
\end{figure}

\begin{figure}[!t]
	\centering
	\begin{subfigure}
		\centering
		\includegraphics[scale=0.52]{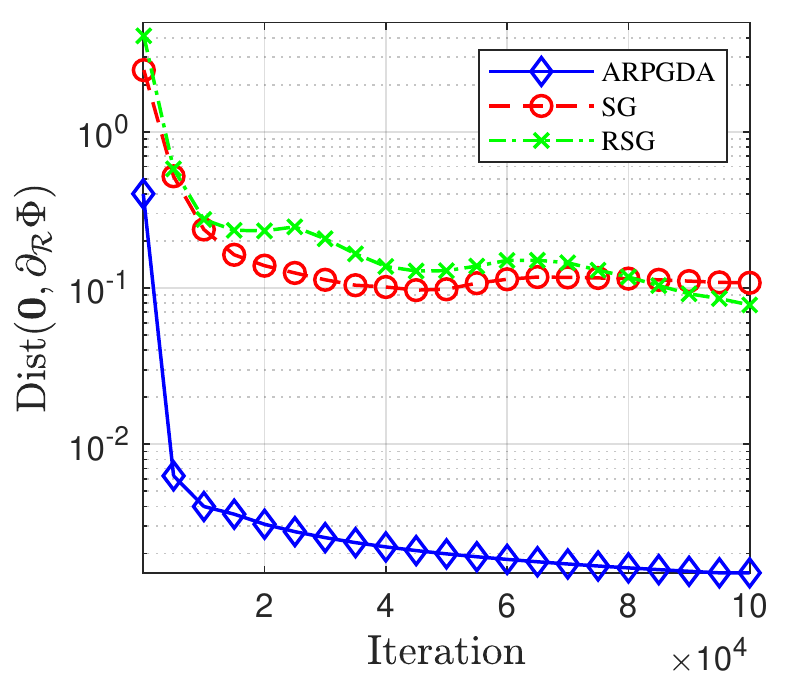}
		\vspace{-7pt}
	\end{subfigure}
	\begin{subfigure}
		\centering
		\includegraphics[scale=0.52]{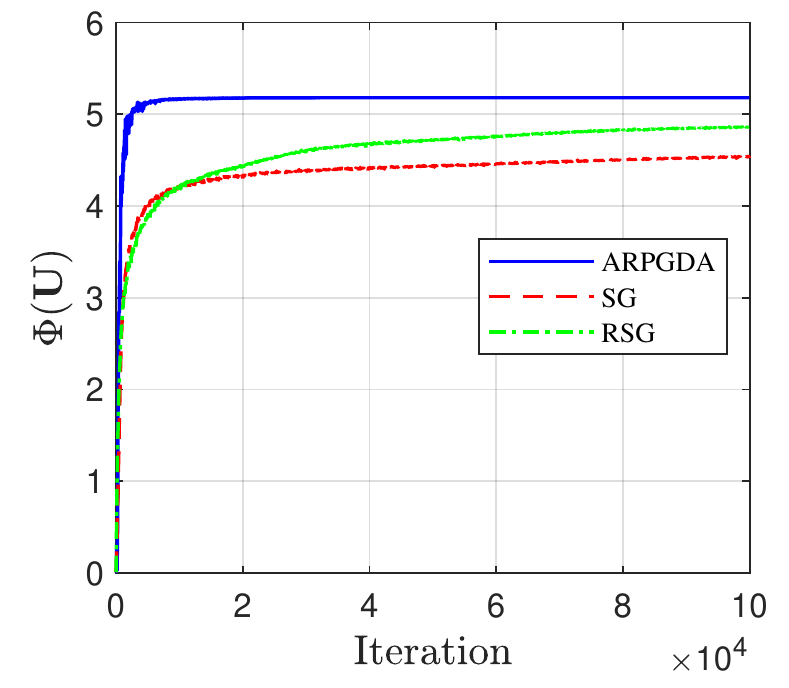}
		\vspace{-7pt}
	\end{subfigure}
	\caption{Performance comparison of an instance on the synthetic data when $ r=2$.  %\rrev{Bo: Some concerns?} 
		%The right figure compares their respective objective function values.
	}
	\vspace{-10pt}
	\label{Synpic2}
\end{figure}

%with zero mean and unit variance \rev{\cite{zalcberg2021fair}}.  
%to control the dimensions and the scale of the problem. 
%\rev{The parameters used in our proposed ARPGDA algorithm are as follows.} %, and show the objective function value of one of the tests to illustrate the performance.
%As in \rrev{\cite{pan2021efficient} Bo: should it be this reference???}, we use a more practical stopping rule for ARPGDA, i.e., $ \|\mathrm{grad}_{\bU} f(\bU,\by)
%\|_{\bU} \leq\varepsilon_1, \|\by - \mathcal{P}_\mathcal{Y}\left(\by_{}+ \nabla_{\by} f\left(\bU_{}, \by_{}\right)\right) \| \leq\varepsilon_2$, and the stopping rule
%\rev{We stop the subgradient methods when the relative is the minimum norm of the (Riemannian) subgradient is less than  $ \varepsilon_1 $, i.e., $ \mathrm{Dist}(\textbf{0}, \rev{\partial_{\mathcal{R}} \Phi}):= \min_{\by\in\Delta} \|\sum_{i\in\mathcal{A}}y_i\mathrm{grad}f_i\|\leq\varepsilon_1$, where $ \mathcal{A}:=\{i\ |\ f_i-\min_i f_i \leq0.1\min_i f_i, i=1,\ldots,n \} $.}
%, is less than  $ \varepsilon_3 $.
%We set \rrev{$ \varepsilon_1 = 0.1$???} and set $ \varepsilon_2 $ small enough such that the solutions of ARPGDA are no worse than those of the subgradient methods.
%\rev{We compare the average objective function value and the average computational time over 10 runs for target dimensions $ r = 1,\ldots,6 $.}
The results on synthetic datasets, wherein the samples $ \{\bx_i\}_{i=1}^n$ with $d=n=200$ are independently generated according to the standard Gaussian distribution \cite{zalcberg2021fair}, are plotted in Fig. \ref{Synpic1}, where $ \hat{\Phi} $ denotes the maximum objective value obtained by the three algorithms.  
From the figure, we can see that i) our proposed ARPGDA is not only much faster than RSG and SG but also can return much better solutions; ii) RSG can return better solutions than SG.  More interestingly, the speed advantage of ARPGDA becomes larger when $r$ increases.
In our tests, we also observe that both of the subgradient methods often reach the preset maximum number of iterations. To further understand the convergence behavior of the three algorithms, we take an instance with $r = 2$ and run all three algorithms till $10^5$ iterations.  The results are plotted in Fig. \ref{Synpic2}. 
%\rev{We stop the subgradient methods when the relative is the minimum norm of the (Riemannian) subgradient is less than  $ \varepsilon_1 $, i.e., $ \mathrm{Dist}(\textbf{0}, \rev{\partial_{\mathcal{R}} \Phi}):= \min_{\by\in\Delta} \|\sum_{i\in\mathcal{A}}y_i\mathrm{grad}f_i\|\leq\varepsilon_1$, where $ \mathcal{A}:=\{i\ |\ f_i-\min_i f_i \leq0.1\min_i f_i, i=1,\ldots,n \} $.}
From this figure, we see that ARPGDA can reduce the approximate minimum norm of the (Riemannian) subgradient, denoted by  $\mathrm{Dist}(\textbf{0}, \partial_{\mathcal{R}}\Phi)$, and improve the objective value much faster than the subgradient-type methods. Here, $\mathrm{Dist}(\textbf{0}, \partial_{\mathcal{R}}\Phi)\coloneqq\min_{\sum_{i \in \mathcal{A}} y_i = 1, y_i \geq 0} \|\sum_{i\in\mathcal{A}}y_i\mathrm{grad}f_i\|$, where $ \mathcal{A}:=\{i \in \{1, 2, \ldots, n\}\ |\ f_i-\min_i f_i \leq0.1\min_i f_i \}.$ %\rrev{Bo: Shall we not compare this term?}


%The results are presented in Fig.\ref{Synpic1}. Remarkably, the two subgradient methods share a similar CPU time, while our ARPGDA algorithm is more efficient to get a better-quality solution. 
%The right figure in Fig. \ref{Synpic1} illustrates our ARPGDA can return a better-quality soltion. 
%We also note that the subgradient methods reach the maximum number of iterations frequently. More specifically,  Fig.\ref{Synpic2} exploits that ARPGDA needs less iterations to find a high-quality solution. We can further analyze that the objective function is usually non-smooth at the stationary points of the FPCA problems. As a consequence, the subgradient methods have to suffer low efficiency near the stationary points. 
%Additionally, the computational complexity increases with the target rank $ r $, for the optimization variable $ \bU $ is a $ d\times r $ matrix.  

\begin{figure}[t]
	\centering
	\begin{subfigure}
		\centering
		\includegraphics[scale=0.52]{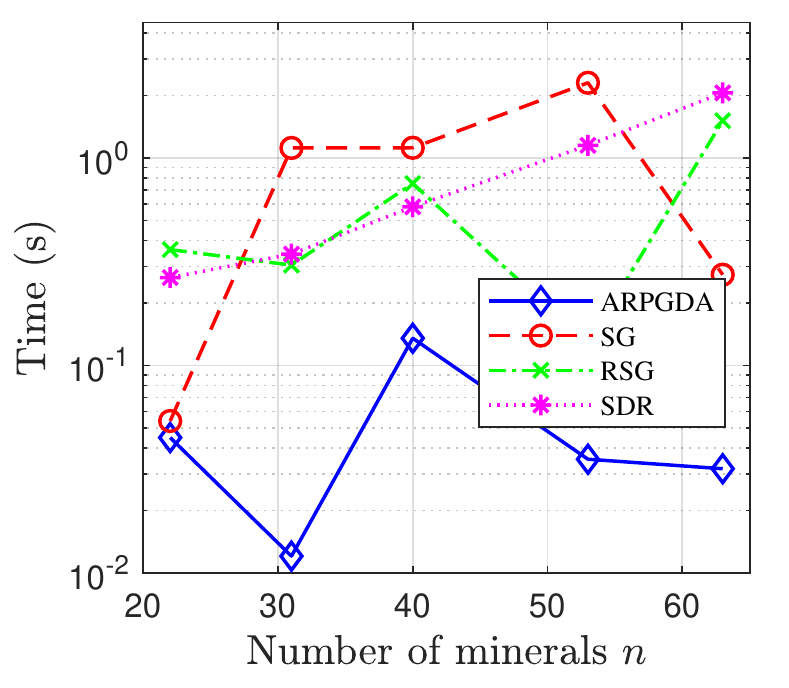}
		\vspace{-7pt}
	\end{subfigure}
	\begin{subfigure}
		\centering
		\includegraphics[scale=0.52]{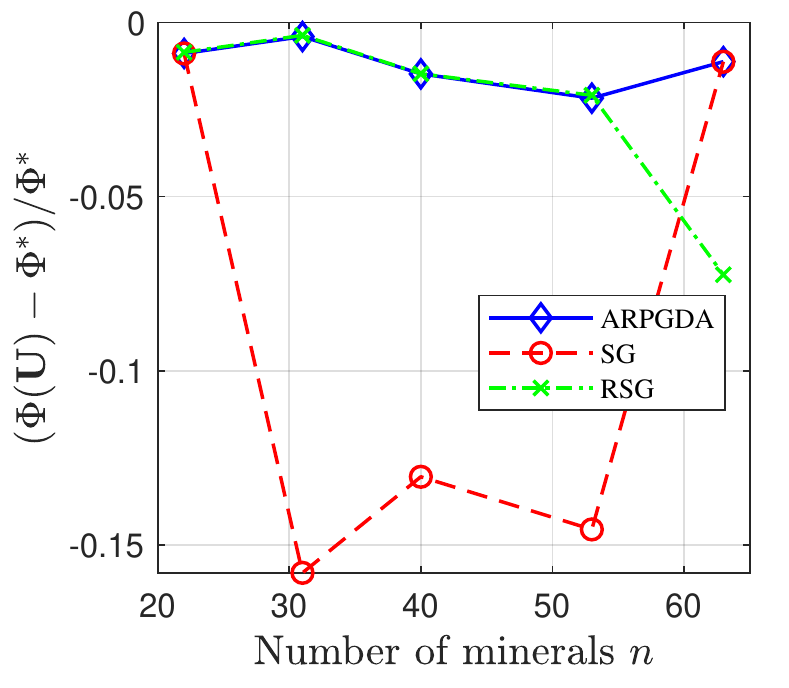}
		\vspace{-7pt}
	\end{subfigure}
	\caption{Average performance comparison on the mineral dataset.}
	\vspace{-7pt}
	\label{Mpic}
\end{figure}

To get more convincing and realistic results,  we consider the spectral signatures of minerals from the
%https://crustal.usgs.gov/speclab/QueryAll07a.php?page=1
 %\href{https://crustal.usgs.gov/speclab/QueryAll07a.php?page=1}{https://crustal.usgs.gov/speclab/QueryAll07a.php?page=1} 
Spectral Library of United States Geological Survey (USGS) 
as in \cite{zalcberg2021fair}. We select a total of $n = 63$ kinds of minerals, each with a reflectance signal $ \{\bx_i\}_{i=1}^{63}$ with $d = 421$ after dropping the signal vectors, whose norms are so small that they may cause a trivial degenerate problem. The vectors are normalized and centered. %We set $\rrev{\varepsilon=0.01????}, \rho=1.1 $ and $\mu_k=5n^2$ \rev{in our proposed ARPGDA algorithm}. 
%We choose the  stopping parameter $\varepsilon_1 = 0.05$ for (R)SG and ARPGDA.
% We still compare the average CPU time over 10 iterations, and set the  stopping parameter $\epsilon_1 = 0.05$ for (R)SG and ARPGDA.
We also present the performance of the SDR approach, where the corresponding semidefinite program is solved via CVX \cite{grant2014cvx} as in \cite{zalcberg2021fair}. 
%since the dimension is not too high. 
Let $\Phi^*$ denote the optimal value of the SDR. By \cite[Theorem 5.2]{samadi2018price}, we know that  this value is equal to the optimal value of the FPCA problem \eqref{FPCAblock}. The averaged results are presented in Fig. \ref{Mpic}. Note that
our ARPGDA algorithm is the fastest among the four algorithms. In particular, except for $n=20$ and $n=40$, ARPGDA is generally 10 times faster than the subgradient-type methods and the SDR approach.  Besides, the objective values returned by ARPGDA are still much better than those returned by SG and RSG on average. 
%It should be mentioned that for $n = 31$ and $n =  63$ our proposed ARPGDA algorithm can find the global optimal solutions. 
We also observe that the local maxima vary dramatically for this dataset, which causes the choppy  curves.
%As shown in Fig.\ref{Mpic}, our proposed ARPGDA algorithm is still efficient for realistic mineral dataset.
%We also note that the quality of the solutions obtained by these algorithms (excepts SDR) is consistently high, which implies that the landscape of this realistic ptoblem is benign.
%We can conclude although ARPGDA has a bit more computational complexity than the two subgradient methods for a single iteration, it has advantage on the total cost of time. 
%However, the SDR approach suffers a much Tr time complexity, needless to say the challenge to get a feasible projection. 
%It is also notable that the time complexity of ARPGDA will not increase rapidly with the number of minerals, meaning that it is suitable for T-scale data.

%\begin{figure}[t]
%	\centering
%	\includegraphics[scale=0.55]{Mtimepic2.pdf}
%	\vspace{-10pt}
%	\caption{The time cost comparison on the mineral dataset.}
%	\label{Mtimepic}
%\end{figure} 
%\begin{figure}[t]
%	\centering
%	\includegraphics[scale=0.64]{Creditpic.pdf}
%	\vspace{-10pt}
%	\caption{The marginal loss of different target ranks $r=1,2,\ldots,20$.}
%	\label{Creditpic}
%\end{figure} 
\begin{figure}[!t]
	\centering
	\begin{subfigure}
		\centering
		\includegraphics[scale=0.52]{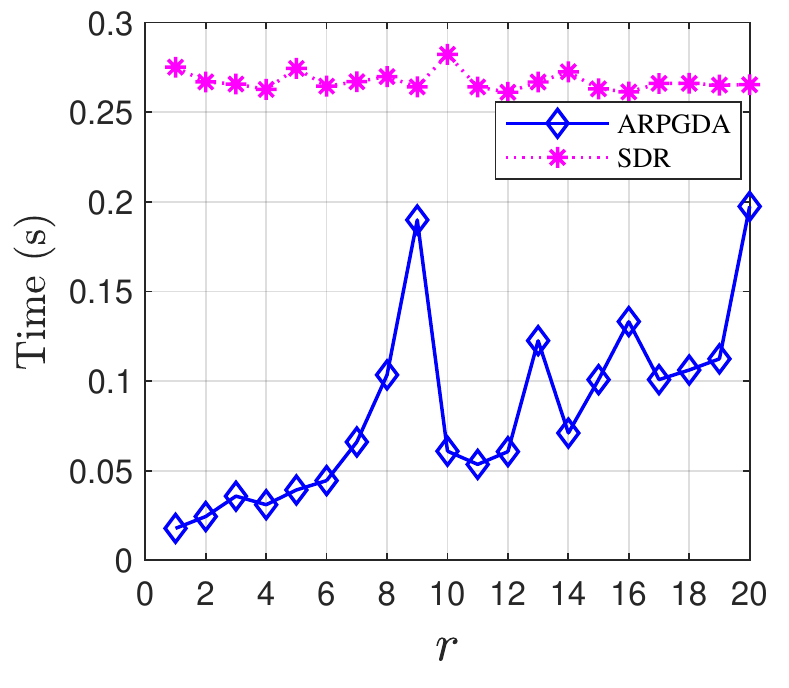}
		\vspace{-7pt}
	\end{subfigure}
	\begin{subfigure}
		\centering
		\includegraphics[scale=0.52]{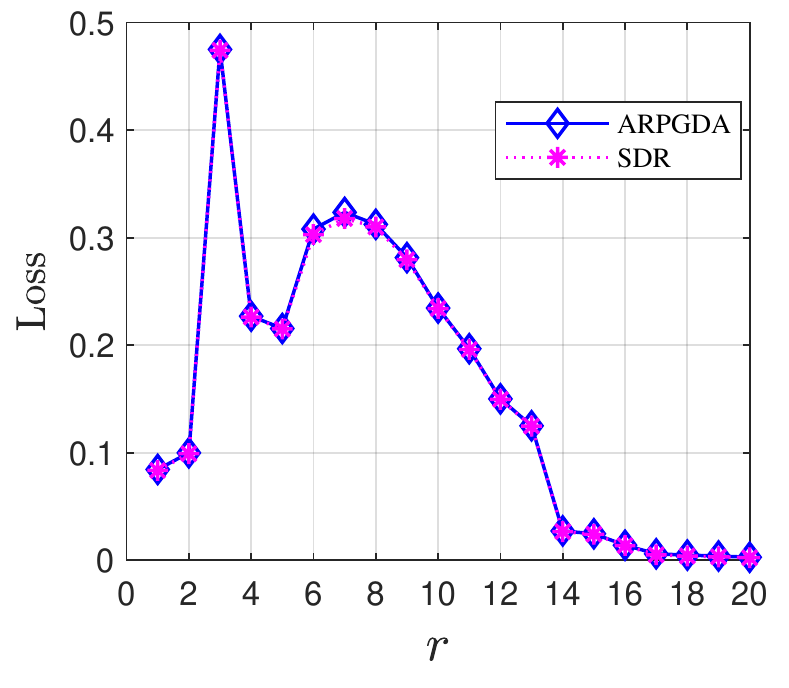}
		\vspace{-7pt}
	\end{subfigure}
	\caption{Average performance comparison on the credit dataset.}
	\vspace{-10pt}
	\label{Creditpic}
\end{figure} 
%Figs. \ref{Synpic} and \ref{Mtimepic} illustrate that our proposed ARPGDA algorithm is more efficient on both the synthetic data and the real world data, and SG and RSG share a simliar performance. More specifically, as shown in Fig. \ref{Synpic}, to get a high-quality subspace $ U $, the two subgradient methods need remakably more iterations than our proposed ARPGDA algorithm. From Fig. \ref{Mtimepic}, we can conclude despite the fact that our ARPGDA algorithm has a little more computational complexity than the two subgradient methods for a single iteration, it still has advantage in the total cost of time. However, the SDR approach via CVXPY \cite{diamond2016cvxpy} suffers a much Tr time complexity due to its complicated matrix variable, needless to say the solution to SDR might not be a feasible projection. It is also notable that as the number of minerals increases, the time complexity of ARPGDA will not increase rapidly, meaning that it is suitable for T-scale data.  

\textbf{Numerical Results on Block FPCA Problem \eqref{FPCAblock} with $n \ll N$.}
The block FPCA problem (\ref{FPCAblock}) with $n \ll N$ arises in a fair machine learning application \cite{samadi2018price,tantipongpipat2019multi}. We choose the Default Credit dataset \cite{yeh2009comparisons}, which contains $ N=30000$ client samples, each with $d = 23$ features including payments and other personal details. The dataset is partitioned into $ n=4 $ sensitive groups according to education and gender. Since the dimension $ d $ is not high and the solutions of the SDR approach are shown to be optimal at all ranks except $ r=8 $ \cite{tantipongpipat2019multi}, the SDR approach is more competitive in this case than the subgradient-type methods. Therefore, we only compare the performance of ARPGDA and the SDR approach. 
The parameters of ARPGDA are $\varepsilon = 10^{-3}$, $\rho = 1.01$, $\theta=1.99$, and $\mu_k \equiv 200 n^2\sqrt{r}$.
%\rev{It is worth mentioning that we can always let $ c_k=0 $ to improve the performance of ARPGDA in our tests, and $ c_k $ is only necessary for some extreme cases.?????}
%Note that for this dataset, choosing $ c_k = 0 $ gives a better performance.
We present  
%which denotes the maximum reconstruction error among the four groups, 
the average computational time and marginal loss value defined in \cite{tantipongpipat2019multi} of SDR and ARPGDA in Fig. \ref{Creditpic}.  
%The bottom figure of Fig.\ref{Creditpic} illustrates that our ARPGDA algorithm is still more efficient, especially when $ r $ is small. 
The left subfigure of Fig. \ref{Creditpic} illustrates that  ARPGDA is always faster than SDR. More specifically, it is always more than 2 times faster than SDR except for the cases when $r = 9, 20$.  For some cases when $r \leq 6$ 
%and $10 \leq r \leq 14$
, it is even more than 5 times faster than SDR. 
As shown in the right subfigure of Fig. \ref{Creditpic}, ARPGDA reaches the SDR upper bound with high accuracy for most of the rank cases, indicating that ARPGDA is effective for block FPCA as well. %We also note that SDR can return global solutions when $ r=6,7 $, while ARPGDA fails, just returning local solutions.

At last, we can conclude from the simulation results that our proposed ARPGDA algorithm is quite efficient in solving various FPCA problems and exhibits a good capability of finding high-quality solutions compared to existing state-of-the-art algorithms \cite{tantipongpipat2019multi,zalcberg2021fair,li2021weakly}.
%At the same time, we have to mention that there is no guarantee that the ARPGDA algorithm can reach a global optimum, despite the fact that we find high-quality solutions in almost all of our experiments.  

\section{Conclusion}
In this paper, we proposed a novel alternating Riemannian/projected gradient descent ascent algorithm for solving the FPCA problem and the more general non-convex linear minimax problem over a Riemannian manifold. We proved that the proposed algorithm can find an $\varepsilon$-stationary point of the above problems within $\mathcal{O}(\varepsilon^{-3})$ iterations. 
We showed via simulation results that the proposed algorithm has a better performance in terms of solution quality and speed when solving the FPCA problems arising in signal processing and machine learning applications.
%We showed via simulation results that the proposed algorithm is quite efficient in solving various FPCA problems arising in signal processing and machine learning and exhibits a good capability of finding high-quality solutions compared to existing state-of-the-art algorithms.
%\bibliographystyle{IEEEbib}
\bibliographystyle{IEEEtran}
%\bibliography{chenbib1007}
\bibliography{reference}
\end{document}